\documentclass[12pt]{article}
\usepackage{amsmath}
\usepackage{graphicx}
\usepackage{enumerate}
\usepackage{natbib}
\usepackage{url} % not crucial - just used below for the URL

\newcommand{\blind}{1}

\def\T{{ \mathrm{\scriptscriptstyle T} }}
\newtheorem{theorem}{Theorem}

\newtheorem{remark}{Remark}

\begin{document}

\def\spacingset#1{\renewcommand{\baselinestretch}%
{#1}\small\normalsize} \spacingset{1}

%%%%%%%%%%%%%%%%%%%%%%%%%%%%%%%%%%%%%%%%%%%%%%%%%%%%%%%%%%%%%%%%%%%%%%%%%%%%%%

\if1\blind
{
  \title{\bf Gaussian variational approximation with composite likelihood for crossed random effect models}
  \author{Libai Xu\thanks{}\hspace{.2cm}\\
    School of Mathematical Sciences, Soochow University, Suzhou, Jiangsu, China,\\
1 Shizi Street, Suzhou, Jiangsu, 215006, China.\\
lbxu@suda.edu.cn\\
    Nancy Reid and Dehan Kong \\
    Department of Statistical Sciences, University of Toronto, Toronto, Canada,\\
    700 University Ave, 9th Floor, Toronto, Ontario M5G 1Z5, Canada.\\
    nancym.reid@utoronto.ca and dehan.kong@utoronto.ca}
  \maketitle
} \fi

\if0\blind
{
  \bigskip
  \bigskip
  \bigskip
  \begin{center}
    {\LARGE\bf Title}
\end{center}
  \medskip
} \fi

\bigskip
\begin{abstract}
Composite likelihood usually ignores dependencies among response components, while
variational approximation to likelihood ignores
dependencies among parameter components.
We derive a Gaussian variational
approximation to the composite log-likelihood function for Poisson  and Gamma regression models with crossed random effects. We show consistency and asymptotic normality of the estimates derived from this approximation and support this theory with some simulation studies.
The approach is computationally much faster than a Gaussian variational approximation to the full log-likelihood function.
\end{abstract}

\noindent%
\vfill

\newpage
\spacingset{1.9} % DON'T change the spacing!
\section{Introduction}
\label{sec:intro}

Generalized linear mixed models with crossed random effects %are a powerful class of models that relate a response variable to categorical predictors by a link function. They
are useful for the analysis and inference
of cross-classified data, such as arise in  educational studies \citep{chung2021cross,menictas2022streamlined}, psychometric research studies
\citep{baayen2008mixed,jeon2017variational}, and medical studies \citep{coull2001crossed}, among others.

%\nancy{
Suppose $Y_{ij}, i = 1, \dots, m; j = 1, \dots, n$ are conditionally independent, given random effects $U_i$ and $V_j$. A generalized linear model for $Y_{ij}$ has density function
\begin{equation}
f(Y_{ij}\mid X_{ij},U_{i},V_{j})=\exp\left[\{Y_{ij}\theta_{ij}-b(\theta_{ij})\}/a(\phi)+c(Y_{ij},\phi)\right], \label{n1}
\end{equation}
and we relate this to the covariates in the usual manner:
\begin{equation}
g(\mu_{ij})= X^{\T}_{ij}\beta + U_{i} + V_{j}, \label{m2}
\end{equation}
where $\mu_{ij} = E(Y_{ij})$.
In \eqref{m2} $X_{ij}=(1,X_{ij})^{\T}$, $\beta$ is a $(p+1)$-vector of fixed-effects parameters, and $U_{i}$ and $V_j$ are independent random effects assumed to follow $N(0,\sigma^{2}_{u})$ and
$N(0,\sigma^{2}_{v})$ distributions, respectively.

%\nancy{
In this paper, we develop a Gaussian variational approximation to a form of composite likelihood for model (\ref{n1}), in order to make the computations more tractable than in a full likelihood approach. We focus on two examples: %{\bf
Poisson regression
%}
, where $Y_{ij} \sim \text{Poisson}(\mu_{ij})$, with $\theta_{ij} = \log(\mu_{ij})$, $b(\theta_{ij}) = \exp(\theta_{ij})$,  $g(\mu_{ij}) = \theta_{ij}$, and $\phi = 1$, and %{\bf
Gamma regression, where $Y_{ij}$ distributed as a Gamma distribution with shape parameter $\alpha = \phi^{-1}$ and expectation parameter $\mu_{ij}$, so that $\theta_{ij} = -\mu_{ij}^{-1}$, $b(\theta_{ij}) =\log(\mu_{ij})$, and
$g(\mu_{ij})=\log(\mu_{ij})$.
%}
%\nancy{

Other approaches to simplifying the likelihood for crossed random effects have been proposed. Penalized quasi-likelihood is discussed in \citet{breslow1993approximate, schall1991estimation}, based on Laplace approximation, although for binary data the resulting estimates are not guaranteed to be consistent \citep{lin1996bias}. \cite{coull2001crossed} developed a Monte Carlo Newton-Raphson algorithm based on  expectation-maximization \citep{mcculloch1997maximum}, and
\cite{ghosh2022backfitting} proposed a backfitting algorithm for  a Gaussian linear model, where the integral can be evaluated explicitly. This approach was extended to logistic regression in \cite{ghosh2022scalable}.

Composite likelihood methods have also been proposed to simplify computation in models with random effects. \cite{renard2004pairwise} and \cite{bartolucci2016pairwise} proposed pairwise likelihood for nested random effects models, and \cite{bellio2005pairwise} developed a pairwise likelihood approach to crossed random effects.

Our approach combines composite likelihood with variational inference. The theory is developed by extending results of \cite{ormerod2012gaussian}, \cite{hall2011theory}, and \cite{hall2011theory} which treat Poisson models with nested random effects. Related variational approaches are developed in \cite{menictas2022streamlined} for linear models, and in \cite{rijmen2013fitting}, \cite{jeon2017variational}, \cite{hui2017variational}, \cite{hui2019semiparametric} for generalized linear mixed models.

We introduce a row-column composite likelihood that requires the computation of only one-dimensional integrals and apply a Gaussian variational approximation to this composite likelihood.
We prove consistency and asymptotic normality of the variational estimates for the Poisson  and Gamma regression models.  Simulations demonstrate that our method is faster than the Gaussian variational approximation to the full likelihood function.

\section{Gaussian variational approximation}

The log-likelihood function for the generalized linear model with crossed random effects is constructed from the marginal distribution of the responses, $Y_{ij}$, given the covariates $X_{ij}$, and depends on the parameters $\Psi = (\beta,\sigma^{2}_{u},\sigma^{2}_{v})$:
\begin{eqnarray}\label{marlikelihood}
\ell(\Psi) &=& \log \int_{R^{m + n}} \prod^{m}_{i=1}\prod^{n}_{j=1}f(Y_{ij}\mid X_{ij},U_{i},V_{j})
p(U_{i})p(V_{j})dU_{i}dV_{j},
\end{eqnarray}
with maximum likelihood estimator
$\hat{\Psi}
=\arg\max_{\Psi}\ell(\Psi) $.
To simplify the computation of the $(m+n)$-dimensional integral
we apply a Gaussian variational approximation to \eqref{marlikelihood} by introducing
pairs of variational parameters $(\mathbf{\mu}_{u_{i}},\lambda_{u_{i}})
, i=1,...,m$ and $(\mathbf{\mu}_{v_{j}},\lambda_{v_{j}}),
j=1,...,n$.
By Jensen's inequality and concavity of the logarithm function:
\begin{eqnarray}
\ell(\Psi)&=& \log \int_{R^{m+n}} \prod^{m}_{i=1}\prod^{n}_{j=1}\left\{f(Y_{ij}\mid X_{ij},U_{i},V_{j})
p(U_{i})p(V_{j}) \frac{\phi(U_{i})\phi(V_{j})}{\phi(U_{i})\phi(V_{j})}dU_{i} dV_{j}\right\}\nonumber\\
&\geq& \int_{R^{m+n}}  \sum^{m}_{i=1}\sum^{n}_{j=1}\phi(U_{i})\phi(V_{j}) \log f(Y_{ij}\mid X_{ij},U_{i},V_{j})  dU_{i} dV_{j} \nonumber\\
&&+
\int_{R^{m}}\sum^{m}_{i=1}\phi(U_{i})\log p(U_{i}) dU_{i} + \int_{R^{n}}\sum^{n}_{j=1} \phi(V_{j})\log p(V_{j}) dV_{j}\nonumber\\
&&- \int_{R^{m}}\sum^{m}_{i=1}\phi(U_{i}) \log \phi(U_{i}) dU_{i} - \int_{R^{n}}\sum^{n}_{j=1} \phi(V_{j})\log \phi(V_{j})dV_{j}=\underline{\ell}(\Psi,\xi),\nonumber
\end{eqnarray}
where
$\phi(U_i)$ and $\phi(V_j)$ are Gaussian densities with means $\mu_{u_{i}}$, $\mu_{v_{j}}$, and variances $\lambda_{u_{i}}$, $\lambda_{v_{j}}$, respectively.
The function $\underline{\ell}(\Psi,\xi)$ is the
 variational lower bound on $\ell(\Psi)$; the variational parameters are $\xi=(\mu_{u_{1}},...,\mu_{u_{m}},\mu_{v_{1}},...,\mu_{v_{n}},\lambda_{u_{i}},...,\lambda_{u_{m}},\lambda_{v_{1}},...,\lambda_{v_{n}})^{\T}$.
 Ignoring some constants, this lower bound simplifies to
\begin{eqnarray}
\underline{\ell}(\Psi,\xi)
&=& \sum^{m}_{i=1}\sum^{n}_{j=1}
\left[Y_{ij}E_{\phi(u_{i})\phi(v_{j})}(\theta_{ij})
-E_{\phi(u_{i})\phi(v_{j})}\{b(\theta_{ij})\}
\right]\nonumber\\
&& + (1/2)\sum^{m}_{i=1}\left\{\log (\lambda_{u_{i}}/\sigma^{2}_{u})-(\mu^{2}_{u_{i}}+\lambda_{u_{i}})/\sigma^{2}_{u}\right\}\nonumber\\
&&+ (1/2)\sum^{n}_{j=1}\{\log (\lambda_{v_{j}}/\sigma^{2}_{v})- (\mu^{2}_{v_{j}}+\lambda_{v_{j}})/\sigma^{2}_{v}\}.
\label{gval}
\end{eqnarray}
The Gaussian variational approximation estimators
maximize the parameters of the  evidence lower bound;
$(\widehat{\underline{\Psi}},\widehat{\underline{\xi}})={\arg\max}_{\Psi,\xi}\underline{\ell}(\Psi,\xi).$
The advantage of using the variational lower bound $\underline{\ell}(\Psi,\xi)$ over $\ell(\Psi)$ is that the former only contains terms
$E_{\phi(u_{i})\phi(v_{j})}(\theta_{ij})$ or
$E_{\phi(u_{i})\phi(v_{j})}\{b(\theta_{ij})\}$ involving a two-dimensional integral, and in our models these can be evaluated explicitly.
For the Poisson regression model,
\begin{eqnarray}
E_{\phi(u_{i})\phi(v_{j})}(\theta_{ij})&=&X^{\T}_{ij}\beta+\mu_{u_{i}}+\mu_{v_{j}},\nonumber\\
E_{\phi(u_{i})\phi(v_{j})}\{b(\theta_{ij})\}&=&
\exp(X^{\T}_{ij}\beta+\mu_{u_{i}}+\lambda_{u_{i}}/2 + \mu_{v_{j}}+\lambda_{v_{j}}/2).\nonumber
\end{eqnarray}
For the Gamma regression model,
\begin{eqnarray}
E_{\phi(u_{i})\phi(v_{j})}(\theta_{ij})&=&-\exp(-X^{\T}_{ij}\beta-\mu_{u_{i}}+\lambda_{u_{i}}/2-\mu_{v_{j}}+\lambda_{v_{j}}/2),\nonumber\\ E_{\phi(u_{i})\phi(v_{j})}\{b(\theta_{ij})\}&=&X^{\T}_{ij}\beta+\mu_{u_{i}}+\mu_{v_{j}}.\nonumber
\end{eqnarray}
%The availability of these explicit expressions  simplified the calculation of $\underline{\ell}(\Psi,\xi)$, compared to $\ell(\Psi)$.

\section{Composite likelihood Gaussian variational approximation}
%\kong{May first need to explain why the GVA of likelihood is %computational expensive. The current logic does not follow smoothly as you were mentioning the advantage of GVA in the last paragraph.}

%Though the GVA of likelihood can simply the calculation, and
%the computational time of the GVA to likelihood increases linearly with the %sample size from the simluation in \cite{ormerod2012gaussian} and types of %random effects, but as the number of random effects increases, GVA will lose %its computational advantage.

%The Gaussian variational approximation to the likelihood function for crossed random effects models can be programmed, but have the additional burden of calculating $mn$ two-dimensional integrals each iteration. Even if there are explicit expressions of these integrals in Poisson and Gamma regression cases, $mn$ times of calculation is needed.  %{\color{blue}{Addressed.}}
In general models, computation of the variational approximation requires $mn$ two-dimensional integrals; for the Poisson and Gamma regression models this simplifies to $mn$ function evaluations. To reduce computation further we develop a version of composite likelihood by considering the two random effects $U_i$ and $V_j$ separately.

First, we ignore the dependence among rows, $Y^{(1)}_{i}=(Y_{i1},...,Y_{in})
, i=1,..., m$, by ignoring the column random effects $V_{j}$'s; i.e., in  \eqref{m2} we replace $X^{\T}_{ij}\beta +U_{i}+V_{j}$ with $X^{\T}_{ij}\beta^r +U_{i}$ , and define the row-composite likelihood function by
\begin{eqnarray}\label{cl1}
CL_{1}(\beta^r,\sigma^{2}_{u})
%&=&\prod^{m}_{i=1}f(Y_{i1},...,Y_{in})\nonumber\\
%&=&\prod^{m}_{i=1}\int_{\mathbf{R}} %f(Y_{i1},...,Y_{in}|X_{ij,U_{i})p(U_{i})dU_{i}\n%onumber\\
&=&\prod^{m}_{i=1}\int_{R} \prod^{n}_{j=1}f(Y_{ij}\mid X_{ij},U_{i})p(U_{i})dU_{i}.\nonumber
\end{eqnarray}
The row-composite log-likelihood function is
\begin{eqnarray}\label{rowc}
c\ell_{1}(\beta^{r},\sigma^{2}_{u})%&=&c\ell_{1}(\tilde{\beta}^{O_{1}},\beta^{O_{2}},\sigma^{2}_{u})
%&=& \sum^{n}_{j=1}\log p(Y_{1j},...,Y_{mj})
= \sum^{m}_{i=1} \log \int_{R} \prod^{n}_{j=1}p(Y_{ij}\mid X_{ij},U_{i})p(U_{i})dU_{i}.
\end{eqnarray}
By a similar operation ignoring row random effects, we can define the column-composite log-likelihood function by
\begin{eqnarray}\label{columnc}
c\ell_{2}(\beta^{c},\sigma^{2}_{v})%&=&\sum^{m}_{i=1} \log p(Y_{i1},...,Y_{in})
&=& \sum^{n}_{j=1} \log \int_{R} \prod^{m}_{i=1}p(Y_{ij}\mid X_{ij},V_{j})p(V_{j})dV_{j}.
\end{eqnarray}
%\nancy{
The notation $\beta^r=(\beta^{r}_{0},\beta_{1},\dots,\beta_{p})$ and $\beta^c=(\beta^{c}_{0},\beta_{1},\dots,\beta_{p})$ is needed, as the misspecification caused by ignoring the column (or row) random effects changes the intercept from $\beta_0$, say, to $\beta^{r}_{0}=\beta_0 + \sigma^2_v/2$ for the row-composite likelihood and $\beta^{c}_{0}=\beta_0 + \sigma^2_u/2$ for the column-composite likelihood, for both the Poisson and the Gamma regression models. The other components, $\beta_{1},\beta_{2},..,\beta_{p},$ are unchanged. %This is discussed further in \S 4 Remark \ref{remark1}. %}. %[NR? Is this better or is this inaccurate?] } %where $\beta^{c}$ is a substitute parameter for $\beta$.
%where $\beta^{c}=(\tilde{\beta}^{O_{3}},\beta^{O_{4}})$ is correspoding to $\beta=(\beta^{O_{3}},\beta^{O_{4}})$,
%$\tilde{\beta}^{O_{3}}$ is the limit of inconsistent estimator of $\beta^{O_{3}}$,
%$\beta^{O_{4}}$ can be estimated consistently.

Based on \eqref{rowc} and \eqref{columnc}, we propose the misspecified row-column composite log-likelihood function by
\begin{eqnarray}\label{row-col}
c\ell(\Psi^{rc})&=& c\ell_{1}(\beta^{r},\sigma^{2}_{v}) + c\ell_{2}(\beta^{c},\sigma^{2}_{u}),
\end{eqnarray}
where $\Psi^{rc}=(\beta^{r}_{0},\beta^{c}_{0},\beta_{1},\dots,\beta_{p},\sigma^{2}_{u},\sigma^{2}_{v})$.
%$\beta^{r}_{0}$ is the misspecified intercept by ignoring the column random effects, while $\beta^{c}_{0}$ ignores the row random effects, $\beta_{1}$ is the shared slope.
%$\beta^{c}\setminus\beta^{r}$ is the parameter vector in $\beta^{c}$ that is different from $\beta^{r}$.
This definition of misspecified composite likelihood %neglects the column random effects in the row composite likelihoods and the row random effects in the column composite likelihood, respectively. This
reduces the computation of an $(m+n)$-dimensional integral in %the marginal likelihood
\eqref{marlikelihood}
to
$m+n$~one-dimensional integrals.
%which
%is essentially different from the row-column composite likelihood proposed by
%\kong{I did not follow the logic of the following.}
\cite{bartolucci2017composite} proposed a different composite likelihood function in which column (or row) random effects are marginalized over instead.
%We introduce the Gaussian variational approximation to the misspecified row-column composite log-likelihood \eqref{row-col} %to offer a closed-form approximation to $c\ell(\Psi^{rc})$
%in the following subsection.

%{\color{blue}{
%We can remove the following part later. I remember you have a question, can we infer the parameters from $c\ell(\Psi^{rc})$? Let's look at the Poisson case,
%\begin{eqnarray}
%c\ell_{1}&=&\sum^{m}_{i=1} \log \int_{R}\prod^{n}_{j=1}\frac{\lambda^{Y_{ij}}}{Y_{ij}!}e^{-\lambda}dU_{i}\nonumber\\
%&=&\sum^{m}_{i=1} \log \int_{R}\prod^{n}_{j=1}\frac{e^{Y_{ij}(\beta^{r}_{0}+\beta_{1}X_{ij}+U_{i})}}{Y_{ij}!}e^{-e^{\beta_{0}+\beta_{1}X_{ij}+U_{i}}}dU_{i},\nonumber\\
%\end{eqnarray}
%which involves m one-dimensional integrals. How do we calculate these integrals? It seems to require some approximation tools, PQL, Laplace, and AGHQ. While GVA method appears to be par in terms of speed with PQL, Laplace and AGHQ \citep{ormerod2012gaussian}. }}

%\subsection{Gaussian variational approximation to misspecified composite likelihood}\label{subs3.2}
%\subsection{Gaussian variational inference}
%To offer a closed-form approximation to $c\ell(\Psi^{rc})$, we now apply
We construct a Gaussian variational approximation to the misspecified row-column composite likelihood \eqref{row-col} using the same variational distributions as in the previous section, leading to
%by introducing variational distributions $\phi(U_{i})$ with
%variational parameters $(\mu_{u_{i}},\lambda_{u_{i}})$ and $\phi(V_{j})$ with variational parameters $(\mu_{v_{j}},\lambda_{v_{j}})$,
%where $\lambda_{u_{i}}>0$ and $\lambda_{v_{j}}>0$, we have
\begin{eqnarray}
&&c\ell(\Psi^{rc})\nonumber\\
& = & \sum^{n}_{j=1}\log \int_{R} \prod^{m}_{i=1} p(Y_{ij}\mid X_{ij},U_{i})p(U_{i})\phi(U_{i})/\phi(U_{i})dU_{i}\nonumber\\
&&+ \sum^{m}_{i=1}\log \int_{R} \prod^{n}_{j=1} p(Y_{ij}\mid X_{ij},V_{j})p(V_{j})\phi(V_{j})/\phi(V_{j})
dV_{j}\nonumber\\
&\geq &\sum^{n}_{j=1}\sum^{m}_{i=1}
\left[Y_{ij}E_{\phi(U_{i})}(\theta^{r}_{ij})
-E_{\phi(U_{i})}\{b(\theta^{r}_{ij})\}
\right]+(1/2)\sum^{m}_{i=1}\left\{
\log(\lambda_{u_{i}}/\sigma^{2}_{u})-(\mu^{2}_{u_{i}}+\lambda_{u_{i}})/\sigma^{2}_{u}\right\}\nonumber\\
&&+\sum^{n}_{j=1}\sum^{m}_{i=1}
\left[Y_{ij}E_{\phi(V_{j})}(\theta^{c}_{ij})
-E_{\phi(V_{j})}\{b(\theta^{c}_{ij})\}\right]+(1/2)\sum^{n}_{j=1}\left\{
\log(\lambda_{v_{j}}/\sigma^{2}_{v})-(\mu^{2}_{v_{j}}+\lambda_{v_{j}})/\sigma^{2}_{v}\right\}\nonumber\\
&=&\underline{c\ell}(\Psi^{rc},\xi).
\label{mvlb}
\end{eqnarray}
%The advantage of $\underline{c\ell}(\Psi^{rc},\xi)$ over $c\ell(\Psi^{rc})$ is that the intractable integrals in the former are transferred to tractable integrals.

For the Poisson regression model,
\begin{eqnarray}
\theta^{r}_{ij} = X^{\T}_{ij}\beta^{r}+U_{i},~E_{\phi(u_{i})}(\theta^{r}_{ij})= X^{\T}_{ij}\beta^{r}+\mu_{u_{i}},
~E_{\phi(u_{i})}\{b(\theta^{r}_{ij})\} = \exp(X^{\T}_{ij}\beta^{r}+\mu_{u_{i}}+\lambda_{u_{i}}/2);\nonumber\\
\theta^{c}_{ij} = X^{\T}_{ij}\beta^{c}+V_{j},~E_{\phi(v_{j})}(\theta^{c}_{ij})= X^{\T}_{ij}\beta^{c}+\mu_{v_{j}}, ~E_{\phi(v_{j})}\{b(\theta^{c}_{ij})\} = \exp(X^{\T}_{ij}\beta^{c}+\mu_{v_{j}}+\lambda_{v_{j}}/2).\nonumber
\end{eqnarray}

For the Gamma regression model,
\begin{eqnarray}
&&\theta^{r}_{ij} = -\exp(-X^{\T}_{ij}\beta^{r}-U_{i}),~
E_{\phi(u_{i})}(\theta^{r}_{ij}) = -\exp(-X^{\T}_{ij}\beta^{r}-\mu_{u_{i}}+\lambda_{u_{i}}/2),\nonumber\\
&& E_{\phi(u_{i})}\{b(\theta^{r}_{ij})\} = X^{\T}_{ij}\beta^{r}+\mu_{u_{i}};\nonumber\\
&&\theta^{c}_{ij} = -\exp(-X^{\T}_{ij}\beta^{c}-V_{j}),~
E_{\phi(v_{j})}(\theta^{c}_{ij}) = -\exp(-X^{\T}_{ij}\beta^{c}-\mu_{v_{j}}+\lambda_{v_{j}}/2),\nonumber\\
&&E_{\phi(v_{j})}\{b(\theta^{c}_{ij})\} = X^{\T}_{ij}\beta^{c}+\mu_{v_{j}}.\nonumber
\end{eqnarray}
We define the estimators based on this Gaussian variational approximation  by
\begin{eqnarray}
(\widehat{\underline{\Psi}}^{rc},\widehat{\underline{\xi}})= \underset{\Psi^{rc},\xi}{\arg\max}\underline{c\ell}(\Psi^{rc},\xi).\nonumber
\end{eqnarray}
To get the estimator $\widehat{\underline{\Psi}}$ from $\widehat{\underline{\Psi}}^{rc}$, we only need to convert $\widehat{\underline{\beta}^{r}_{0}}$ and
$\widehat{\underline{\beta}^{c}_{0}}$ to $\widehat{\underline{\beta}}_{0}$.
This is discussed in \S 4 Remark 1.
The advantage of $\underline{c\ell}(\Psi^{rc},\xi)$ over $\underline{\ell}(\Psi,\xi)$ is that the former only involves $m+n$ one-dimensional integrals, and especially, for Poisson and Gamma regression models with explicit expressions, even one-dimensional integrals are no longer involved; only $m+n$ calculations are needed. The misspecified row-column composite likelihood-based variational approximation enhances the efficiency of calculations compared to the variational approximation to the log-likelihood function.

\section{Theoretical Properties}
In this section we present our main results on the consistency and convergence rates of the parameter estimates based on the Gaussian variational approximation introduced above. We denote the true value of the parameter with a superscript $0$. %to composite log-likelihood function for the Poisson and Gamma regression models with crossed random effects, respectively. With a bit abuse of the notation, we add a superscript $0$ to denote the true value of a parameter. We consider the following two cases:
\begin{itemize}
    \item Poisson regression model:
\begin{eqnarray}\label{epmm}
&&Y_{ij}\mid X_{ij},U_{i},V_{j}~\text{follows~independent~Poisson~with~mean} ~\exp(\beta^{0}_{0}+\beta^{0}_{1}X_{ij}+U_{i}+V_{j}),\nonumber\\
&& U_{i} \sim N(0,(\sigma^{2}_{u})^{0}), V_{j} \sim N(0,(\sigma^{2}_{v})^{0}), U_{i}~\text{and}~V_{j}~\text{are independent}.
\end{eqnarray}
\item Gamma regression model:
\begin{eqnarray}\label{epga}
&&Y_{ij} \mid X_{ij},U_{i},V_{j}~\text{independent~Gamma~with~$\alpha$ and mean} ~\exp(\beta^{0}_{0}+\beta^{0}_{1}X_{ij}+U_{i}+V_{j}),\nonumber\\
&& U_{i} \sim N(0,(\sigma^{2}_{u})^{0}), V_{j} \sim N(0,(\sigma^{2}_{v})^{0}), U_{i}~\text{and}~V_{j}~\text{are independent},
\end{eqnarray}
where the shape parameter $\alpha$ is assumed known.
\end{itemize}
%\begin{assumption*}[A9]
%$m/n=\Omega(1)$, which means that $m$ and $n$ have the same order.
%\end{assumption*}

%\kong{Our asymptotic theories are under the asymptotic regime $ m/n=\Omega(1)$. This is different from the one used in \cite{hall2011asymptotic}. Under the asymptotic regime $n/m$ tending to zero, the convergence rate will be dominated by $n^{-1/2}$. ... . Then add why Hall's asymptotic regime can not be applied to our case.} {\color{blue}{Addressed.}}{\color{red} Is it correct to describe it as a different asymptotic regime? Hall only had one random effect.} {\color{blue} A good question! Differet models(nested,crossed) have different requirements on $n/m$.}

%Our asymptotic theories are under the asymptotic regime where $ m/n=\Omega(1)$. This is different from the one used in \cite{hall2011asymptotic}. The asymptotic regime $n/m$ tending to zero in \cite{hall2011asymptotic} is in keeping with the number of groups being large compared within sample size for a mixed model with nested effects, as often arises in practice, it is a sufficient condition but not a necessary condition to ensure the asymptotic theories.
%As the asymptotic manner $n/m$ approaches zero, the rate of convergence is dominated by $n^{-1/2}$, indicating that row random effects can be disregarded, in contrast to the cross random effects we take into consideration. Hence, \cite{hall2011asymptotic} asymptotic regime can not be applied to our case.}}

Under the assumptions  {\rm (A1)-(A9)} in the \S S1 of the Supplementary Material, we have the following two theorems.
\begin{theorem}\label{thm1}
As $m$ and $n$ diverge, such that $m/n=\Omega(1)$, for both Poisson and Gamma regression models,
\begin{eqnarray}
\widehat{\underline{\beta}}_{0}-\beta^{0}_{0} &=& O_{p}(m^{-1/2}+n^{-1/2}),~~\widehat{\underline{\beta}}_{1}-\beta^{0}_{1} =O_{p}(m^{-1/2}+ n^{-1/2}),\nonumber\\
\widehat{\underline{\sigma}}^{2}_{u}-(\sigma^{2}_{u})^{0} &=& O_{p}(m^{-1/2}+n^{-1/2}),~\widehat{\underline{\sigma}}^{2}_{v}-(\sigma^{2}_{v})^{0} = O_{p}(m^{-1/2}+n^{-1/2}).\nonumber
\end{eqnarray}
%as $m$ and $n$ diverge to infinity.
\end{theorem}
\begin{remark}\label{remark1}%{\color{red} I don't think this remark belongs here, because one has to turn to the proof to find $\widehat{\beta_0^r}$ and $\widehat{\beta_0^c}$. The remark might say something like: As noted above, the intercept terms for the row-only and column-only composite likelihood functions are not equal to $\beta_0$. In the proof we derive the probability limits of the variational estimates of $\beta^r_0$ and $\beta^c_0$, and show how to estimate $\beta_0$ from these and estimates of $\sigma^2_u$ and $\sigma^2_v$. [if this is correct]} {\color{blue}Yes! It is exactly correct! }
As noted in \S 3, the intercept terms for the row-only and column-only composite likelihood functions are not equal to $\beta_0$. In the proof we derive the probability limits of the variational estimates of $\beta^r_0$ and $\beta^c_0$, and show how to estimate $\beta_0$ from these and estimates of $\sigma^2_u$ and $\sigma^2_v$.
\end{remark}
%To consider the precise asymptotic distributional behavior of composite likelihood Gaussian variational approximation estimates for the Poisson and Gamma regression models with crossed random effects,
In addition, we have the following asymptotic normality results.
\begin{theorem}
For both Poisson and Gamma regression models, as $m$ and $n$ diverge, such that $m/n=\Omega(1)$, we have
\begin{eqnarray}
\widehat{\underline{\beta}}_{0}-\beta^{0}_{0} &=&
m^{-1/2}1^{\T}_{3}Z_{1}1_{3} + n^{-1/2}1^{\T}_{3}Z_{2}1_{3} + o_{p}(m^{-1/2}+n^{-1/2}),\nonumber
\end{eqnarray}
where $1_{3}=(1,1,1)^{\T}$, $Z_{1}$ and $Z_{2}$ are independent normal distributions with mean zero and covariance
\begin{eqnarray}
\Gamma_{1}
&=& \frac{1}{8}\left(
  \begin{array}{ccc}
    2[\exp\{(\sigma^{2}_{u})^{0}\}-1] & 2(\sigma^{2}_{u})^{0} & -\{(\sigma^{2}_{u})^{0}\}^{2} \\
    2 (\sigma^{2}_{u})^{0} & 2(\sigma^{2}_{u})^{0} & 0 \\
     -\{(\sigma^{2}_{u})^{0}\}^{2} & 0 &     \{(\sigma^{2}_{u})^{0}\}^{2}\\
  \end{array}
\right),\nonumber
\end{eqnarray}
and
\begin{eqnarray}
\Gamma_{2}
&=&\frac{1}{8}\left(
  \begin{array}{ccc}
    2[\exp\{(\sigma^{2}_{v})^{0}\}-1] & 2 (\sigma^{2}_{v})^{0} & -\{(\sigma^{2}_{v})^{0}\}^{2} \\
    2(\sigma^{2}_{v})^{0} & 2(\sigma^{2}_{v})^{0} & 0 \\
     -\{(\sigma^{2}_{v})^{0}\}^{2} & 0 & \{(\sigma^{2}_{v})^{0}\}^{2} \\
  \end{array}
\right)\nonumber
\end{eqnarray}
respectively;
\begin{eqnarray}
\widehat{\underline{\sigma}}^{2}_{u}-(\sigma^{2}_{u})^{0} &=& m^{-1/2}Z_{3} + o_{p}(n^{-1/2}+m^{-1/2}),\nonumber
\end{eqnarray}
where the random variable $Z_{3}$ follows $N(0,2\{(\sigma^{2}_{u})^{0}\}^{2})$;
\begin{eqnarray}
\widehat{\underline{\sigma}}^{2}_{v}-(\sigma^{2}_{v})^{0} &=& n^{-1/2}Z_{4} + o_{p}(m^{-1/2}+n^{-1/2}),\nonumber
\end{eqnarray}
where the random variable $Z_{4}$ follows $N(0,2\{(\sigma^{2}_{v})^{0}\}^{2})$.

For Poisson regression, the slope term satisfies
\begin{eqnarray}
\widehat{\underline{\beta}}_{1}-\beta^{0}_{1
}&=& (mn)^{-1/2}Z_{5} + o_{p}\{m^{-1}+(mn)^{-1/2} + n^{-1}\},\nonumber
\end{eqnarray}
the random variable $Z_{5}$ follows a normal distribution with zero mean and variance
\begin{eqnarray}
&&\exp\{-\beta^{0}_{0}-(\sigma^{2}_{u})^{0}/2-(\sigma^{2}_{v})^{0}/2\}\tau_{1}
+\tau^{2}_{1}\tau_{2} 1^{\T}_{2}\Gamma_{3}1_{2},\nonumber
%\frac{\phi^{2}(\beta^{0}_{1})\phi^{2}_{1}(\beta^{0}_{1})[\exp\{(\sigma^{2}_{v})^{0}\}-1]}{4\{\phi_{2}(\beta^{0}_{1})\phi(\beta^{0}_{1})-\phi^{2}_{1}(\beta^{0}_{1})\}^{2}},\nonumber
\end{eqnarray}
%and $Z_{4}$ follows a normal distribution with zero mean and variance
%\begin{eqnarray}
%\frac{\phi^{2}(\beta^{0}_{1})\phi^{2}_{1}(\beta^{0}_{1})[\exp\{(\sigma^{2}_{u})^{0}\}-1]}{4\{\phi_{2}(\beta^{0}_{1})\phi(\beta^{0}_{1})-\phi^{2}_{1}(\beta^{0}_{1})\}^{2}};\nonumber
%\end{eqnarray}
where $\tau_{1}=\phi(\beta^{0}_{1})/\{\phi_{2}(\beta^{0}_{1})\phi(\beta^{0}_{1})-\phi^{2}_{1}(\beta^{0}_{1})\}$, $\tau_{2}=\phi_{2}(2\beta^{0}_{1})-2\phi_{1}(\beta^{0}_{1})\phi_{1}(2\beta^{0}_{1})/\phi(\beta^{0}_{1}) + \phi_{1}(\beta^{0}_{1})^{2}\phi(2\beta^{0}_{1})/\phi(\beta^{0}_{1})^{2}$, and
\begin{eqnarray}
\Gamma_{3}
&=& \frac{1}{4}\left(
  \begin{array}{cc}
    \exp\{(\sigma^{2}_{u})^{0}\}[\exp\{(\sigma^{2}_{v})^{0}\}-1] & [\exp\{(\sigma^{2}_{u})^{0}\}-1][\exp\{(\sigma^{2}_{v})^{0}\}-1]\\
~[\exp\{(\sigma^{2}_{u})^{0}\}-1][\exp\{(\sigma^{2}_{v})^{0}\}-1]&      \exp\{(\sigma^{2}_{v})^{0}\}[\exp\{(\sigma^{2}_{u})^{0}\}-1]\\
  \end{array}
\right);\nonumber
\end{eqnarray}
For Gamma regression, the slope term satisfies
\begin{eqnarray}
\widehat{\underline{\beta}}_{1}-\beta^{0}_{1}=(mn)^{-1/2}1^{\T}_{2}Z_{6}1_{2}+o_{p}\{m^{-1}+(mn)^{-1/2}+n^{-1}\},\nonumber
\end{eqnarray}
where $Z_{6}$ follows a normal distribution with zero mean and covariance
\begin{eqnarray}
\widetilde{\Sigma}=\frac{1}{4\alpha}\left(
                     \begin{array}{cc}
                      \alpha[\exp\{(\sigma^{2}_{v})^{0}\}-1]+\exp\{(\sigma^{2}_{v})^{0}\} &
                        1 \\
                        1  &  \alpha[\exp\{(\sigma^{2}_{u})^{0}\}-1]+\exp\{(\sigma^{2}_{u})^{0}\}
                         \\
                     \end{array}
                   \right).\nonumber
\end{eqnarray}
\end{theorem}

%{\color{blue}
%We assume above that $m/n=\Omega(1)$, whereas  \citet{hall2011asymptotic} have $n/m \rightarrow 0$.  {\color{red}Seems weird to have $m/n$ and then $n/m$. Is this correct? Maybe say something like: We assume above that $m$ and $n$ diverge at the same rate, whereas Hall et al. assume $n/m \rightarrow 0$, to reflect the usual nested effects setting, with many groups (random effects) and relatively few observations per group.} \cite{hall2011asymptotic}'s asymptotic manner corresponds to a situation in a mixed model with nested random effects, where the number of groups is significantly larger than the sample size. {\color{red} In Hall $m$ is the number of groups, whereas $n$ is the number of observations in each group, so what we would normally call the sample size is $mn$. } While \cite{hall2011asymptotic}'s asymptotic manner is sufficient to ensure the validity of the asymptotic theory, it is not necessary. In the generalized linear mixed model with crossed random effects, if $n/m$ tends to $0$, the convergence rate will be dominated by $n^{-1/2}$, indicating that row random effects can be disregarded. However, in the case we consider, the row random effects should not be ignored.} %Therefore, \cite{hall2011asymptotic}'s asymptotic manner is not applicable to our situation.

\begin{remark}
We assume above that $m$ and $n$ diverge at the same rate, whereas \cite{hall2011asymptotic} assume $n/m \rightarrow 0$, to reflect the usual nested effects setting, with many groups (random effects) and relatively few observations per group. %\cite{hall2011asymptotic}'s asymptotic manner corresponds to a situation in a mixed model with nested random effects, where the number of groups is significantly larger than the sample size. {\color{red}
%In \cite{hall2011asymptotic} $m$ is the number of groups, whereas $n$ is the number of observations in each group, so what we would normally call the sample size is $mn$.
While \cite{hall2011asymptotic}'s asymptotic manner is sufficient to ensure the validity of the asymptotic theory, it is not necessary. In the generalized linear mixed model with crossed random effects, if $n/m$ tends to $0$, the convergence rate will be dominated by $n^{-1/2}$, indicating that row random effects can be disregarded. However, in the case we consider, the row random effects should not be ignored.
\end{remark}

\section{Simulation Studies}

In this section, we perform simulations to evaluate the effectiveness of the proposed approach. We assess the performance of both Poisson and Gamma regression models with crossed random effects as specified in equations \eqref{epmm} and \eqref{epga}. For both models, we independently generate $X_{ij}$ from $N(1,1)$ for $ i=1, \ldots, m$ and $ j=1,\ldots, n$. We set $\beta^{0}_{0}=-2$, $\beta^{0}_{1}=-2$, $(\sigma_{u})^{0}=0.5$, and $(\sigma_{v})^{0}=0.5$. For the Gamma regression model, there is an additional shape parameter set as $ \alpha = 0.8$. We consider two different scenarios: $ (m,n)=(50,50)$ and $(m,n)=(100,100)$.

We compare the proposed approach with the Gaussian variational inference based on the full log-likelihood function. We report the mean and standard deviation of the parameter estimates based on 1000 repetitions of Monte Carlo studies. We also calculate the average computational times of both methods, using R version 4.1.1 on a laptop equipped with a 2.8 gigahertz Intel Core i7-1165G7 processor and 16 gigabytes of random access memory. The results are presented in Table \ref{table1}.

{\scriptsize
\begin{table}
\def~{\hphantom{0}}
\caption{ Simulations for Poisson and Gamma regression models with crossed random effects were conducted using parameters $\beta^{0}_{0}=-2$, $\beta^{0}_{1}=-2$, $(\sigma_{u})^{0}=0.5$, and $(\sigma_{v})^{0}=0.5$ under two different scenarios: $(m,n)=(50,50)$ and $(m,n)=(100,100)$. For the Gamma regression model, the additional shape parameter was set as $ \alpha = 0.8$.
``GVA'' denotes the Gaussian variational approximation
 to the log-likelihood function and ``GVACL'' denotes our proposed method. ``Mean(SD)'' reports the mean and standard deviations of the parameter estimates based on the $1000$ datasets; ``MESE'' denotes the average of the standard errors of the estimated parameters calculated using the asymptotic variance formula;
``Mean Time(s)'' is the average time of each simulation in seconds. For the GVA method, we did not report the MESE since no asymptotic variance was derived in previous literature, denoted by "NA". }
\label{table1}\begin{center}
	\begin{tabular}{ccccc}
		\hline	
        Poisson &\multicolumn{2}{c}{$(m,n)=(50,50)$}&\multicolumn{2}{c}{$(m,n)=(100,100)$}\\		\cline{2-3~}\cline{4-5~}
		 & GVA & GVACL & GVA & GVACL\\
         Mean(SD) & -1.99(0.07) & -2.04(0.14) & -1.99(0.08) & -2.03(0.09)\\
         MESE & NA & 0.11 & NA & 0.07\\
         Mean(SD) & -1.99(0.07) & -1.99(0.08) & -1.99(0.03) & -2.00(0.04)\\
          MESE & NA & 0.10 & NA & 0.05\\
         Mean(SD) & 0.47(0.10) & 0.53(0.09) & 0.49(0.05) & 0.52(0.05)\\
          MESE & NA & 0.05 & NA & 0.04\\
         Mean(SD) & 0.46(0.10) & 0.52(0.09) & 0.49(0.05) & 0.52(0.05)\\
          MESE & NA & 0.05 & NA & 0.04\\
        Mean Time(s)& 3.55 & 0.21 & 14.04 & 0.63 \\
        \hline
        Gamma &\multicolumn{2}{c}{$(m,n)=(50,50)$}&\multicolumn{2}{c}{$(m,n)=(100,100)$}\\
        \cline{2-3~}\cline{4-5~}
         Mean(SD) & -2.00(0.10) & -2.01(0.11) & -2.00(0.07) & -2.00(0.07)\\
          MESE & NA & 0.10 & NA & 0.07\\
         Mean(SD) & -2.00(0.02) & -2.00(0.03) & -2.00(0.01) & -2.00(0.01)\\
        MESE & NA & 0.03 & NA & 0.01\\
         Mean(SD) & 0.49(0.05) & 0.50(0.05) & 0.50(0.04) & 0.50(0.04)\\
          MESE & NA & 0.05 & NA & 0.04\\
        Mean(SD) & 0.49(0.06) & 0.50(0.06) & 0.50(0.04) & 0.50(0.04)\\
         MESE & NA & 0.05 & NA & 0.04\\
         Mean Time(s)& 4.77 & 0.28 & 19.00 & 0.88 \\
         \hline
		\end{tabular}
	\end{center}
\end{table}}

Our results demonstrate that both methods yield similar parameter estimates. However, our proposed approach exhibits a notable advantage in terms of computational efficiency compared to the Gaussian variational approximation to the log-likelihood function. We further investigate additional settings to assess the performance of both methods. The results are included in the supplementary \S S6 and Table S1, and the findings are similar.

%will achieve more precision of parameter estimates as $(m,n)$ become larger, and our method loses a little bit of accuracy but
%is computationally more efficient than the Gaussian variational approximation.

\section{Application}
%The data description from \cite{adam2021modeling}, ``
%The motor vehicle insurance of a general insurance %company in Indonesia in 2014 data sourced from
%the Financial Services Authority (FSA) is used. The %amount of the claim is a response variable and the
%deductible is the explanatory variable. Selection of %response variables and explanatory variables refers
%to Pratama [8]. The actual area code is 35 areas, but %only 17 areas are included in the model because in
%the other 18 areas there are no claim events. The area %code is considered as a random effect and the
%month of occurrence is a random effect which is %considered to follow a first-order autoregressive
%process. For the purposes of this research, the data was %taken in part through simple random sampling
%from the actual data of 175,000 items. In each area and %month of the occurrence 10 policy numbers were
%taken as samples.''

We illustrate our methods with data concerning motor vehicle insurance in Indonesia, in 2014, obtained from the Financial Services Authority \citep{adam2021modeling}. The dataset consists of 175,381 claim events from 35 distinct areas over a period of 12 months. Each claim event includes information on the claim amount and deductible. However, 17 out of the 35 areas did not have any claim events during the entire 12-month period. Consequently, our analysis focuses on the remaining 18 areas. Table S3 includes the counts of claim events for each area and month.

We randomly selected one claim event from each area and month. The response variable, denoted as $Y_{ij}$, represents the claim amount, while the explanatory variable, denoted as $X_{ij}$, represents the deductible, where $i=1, \dots, 18$ and $ j=1, \ldots, 12$. A portion of the dataset can be found in Table S4 of the supplementary material. From that table, we can see both $X_{ij}$ and $Y_{ij}$ have large magnitudes, therefore we divide them by $10^7$ to avoid singular Hessian matrices in computation  \citep{adam2021modeling}. %\kong{Be a bit specific about what singular case means.} {\color{blue}{Addressed.}}
We then fit a Gamma regression model, introducing random effects for both area  and month of occurrence.

To evaluate the proposed method and the Gaussian variational approximation to the log-likelihood function, we regenerated $1000$ datasets  using different random seeds at each step of the random event selection process. We then applied both methods and reported the mean and standard deviation of both estimators across these $1000$ datasets. %To validate the asymptotic variance derived in Theorem \ref{thm3}, %\kong{add where, say equation ?? or theorem ??} {\color{blue}{Addressed.}},
%we also report the averages of the estimated standard errors of $\widehat{\beta}_0$ and $\widehat{\beta}_1$.
In addition, we record the average time it takes to fit the model. Table \ref{table3} presents the results of this analysis.
\begin{table}[h]
\caption{Real data results: ``GVA'' is the Gaussian variational approximation
 to the log-likelihood function, and ``GVACL'' is our proposed method. ``Mean Time(s)'' is the average time of each simulation in seconds; ``Mean(SD)'' reports the mean and sample standard deviations of the model parameters based on the $1000$ datasets. %``MESE'' is the average of the standard error of the estimated parameters calculated using the asymptotic variance formula. %For the GVA method, we did not report the MESE, denoted by ``\textsl{NA}'', as no asymptotic variance was derived in previous literature
}\label{table3}\begin{center}
   % \resizebox{\textwidth}{8mm}{
	%		\setlength{\tabcolsep}{2mm}{
				\begin{tabular}{cccccc}
					%\hline
					\hline	&~&$\beta_{0}$&$\beta_{1}$&
     $\sigma_{u}$&
     $\sigma_{v}$\\
 %    \cline{3-4~}\cline{5-6~}\cline{7-8~}\cline{9-10~}
					  Method & Mean Time(s) & Mean(SD) & Mean(SD) & Mean(SD) & Mean(SD)  \\
       \multicolumn{1}{l}{GVA} & 1.19 & -0.71 (0.12) & 3.17 (1.37) & 0.29(0.11)& 0.10(0.10)  \\
%&~&~&\multicolumn{2}{c}{$~$}&\multicolumn{2}{c}{$~$}&~\\
                    \multicolumn{1}{l}{GVACL} & 0.15 & -0.74 (0.12)
                      & 3.60 (1.26)  & 0.29(0.11)  & 0.14(0.12)  \\	
                    \hline
                    %\hline
		\end{tabular}%}}
	\end{center}
\end{table}
From Table \ref{table3}, it is evident that the two methods yield comparable model parameter estimates. However, our proposed approach ``GVACL'' is much faster in terms of computation time.% \kong{The MESE results do not look so good. We may need to discuss whether to include it or not.}

\section{Discussion}\label{discuss}
%In this article, we have extended the Gaussian variational approximation methods to the generalized linear mixed models with crossed random effects and first proposed the composite likelihood Gaussian variational approximation for the generalized linear mixed models with crossed random effects. We derived closed-form approximations fully to the composite likelihood for commonly observed responses,
%inferred the link between variational composite estimates of composite log-likelihood and true parameters, and proved the consistency and asymptotic normality of
%variational estimates for the Poisson response case and Gamma response case. Through simulations, the composite likelihood Gaussian variational approximation is faster than the Gaussian variational approximation and still guarantees the consistency and asymptotic normality of variational estimates for the Poisson mixed model with crossed random effects
%and the Gamma regression model with crossed random effects. Besides,
In this article, we cover two examples: the Poisson regression and Gamma regression.
An interesting future direction is to extend the results to other generalized linear mixed models.  For example, the  logistic regression model with crossed random effects has
\begin{eqnarray}
&&Y_{ij}=1\mid X_{ij},U_{i},V_{j} \sim \text{Bernoulli}[1/\{1+\exp(-\beta^{0}_{0}-\beta^{0}_{1}X_{ij}-U_{i}-V_{j})\}],\nonumber\\ %\label{ebmm}\\
&&U_{i}\sim N(0,(\sigma^{2}_{u})^{0}), V_{j}\sim N(0,(\sigma^{2}_{v})^{0}),U_{i}~\text{and}~V_{j}~\text{are~independent}.\nonumber
\end{eqnarray}
Unlike the Poisson and Gamma regression models, there is no analytic solution to
the two-dimensional integral
\begin{eqnarray}\label{eq2315}
E_{\phi(u_{i})\phi(v_{j})}\{b(\theta_{ij})\}=\int_{R^{2}}\log\{1+\exp(X^{\T}_{ij}\beta+U_{i}+V_{j})\}\phi(U_{i})\phi(V_{j})dU_{i}dV_{j}
\end{eqnarray}
appearing in \eqref{gval}.
%\kong{Add where this term appears in which equation.} {\color{blue}{Addressed}}
%for the logistic regression case.
One solution is to evaluate the integrals using adaptive Gauss-Hermite quadrature \citep{liu1994note} with a product of $N_{1}\times N_{2}$ quadrature points over the two-dimensional integral in \eqref{eq2315}.
%\kong{each two-dimensional reads confusing. Maybe a bit more specific}.
The one-dimensional integrals arising in the row-column composite likelihood
%$E_{\phi(U_{i})}\{b(X^{\T}_{ij}\beta^{r}+U_{i})\}= \int \log\{1+\exp(X^{\T}_{ij}\beta^{r}+U_{i})\}\phi(U_{i})dU_{i}$ and
%$E_{\phi(V_{j})}\{b(X^{\T}_{ij}\beta^{c}+V_{j})\}= \int \log\{1+\exp(X^{\T}_{ij}\beta^{c}+V_{j})\}\phi(V_{j})dV_{j}$ in
%$\underline{c\ell}(\Psi^{rc},\xi)$
of \eqref{mvlb}
%\kong{add where} {\color{blue}{Addressed.}}
 can also be computed by adaptive Gauss-Hermite quadrature. Compared to the Poisson and Gamma cases, the computational time of the proposed method is longer, %because of the integral calculation,
 especially when $m$ and $n$ are large. %Other methods to obtain  a closed form variational approximation to log-likelihood for the logistic regression model, a reparameterization of the models with the help of auxiliary variables \cite{hui2017variational}, .

%One  when $m$ and $n$ are larger, this calculation takes a lot of time.
%In the Bernoulli case, since $E_{p(v_{j})}b(X^{\T}_{ij}\beta+U_{i}+V_{j})$ have no analytic form, we explore
%the $f$ and $g$ through simulations in Section \ref{discuss}.

We carried out some limited simulations for the logistic model
%under setting S-1 %of supplementary Section S6,
%with different intercepts, slopes, and variances.
and report the results in the supplementary \S S6 and Table S2. We found that the computational time of the proposed method is much shorter than for the Gaussian variational approximation to the log-likelihood function. Based on the simulation results, we conjecture that the intercept and slope estimates and the sum of random effect variance estimates have the following relationship:
\begin{eqnarray}
\widehat{\underline{\beta}}_{0}\approx (\widehat{\underline{\beta}^{r}_{0}}+\widehat{\underline{\beta}^{c}_{0}})/[2\{1-0.1(\widehat{\underline{\sigma}}^{2}_{u}+\widehat{\underline{\sigma}}^{2}_{v})\}]
~~\text{and}~~
\widehat{\underline{\beta}}_{1}\approx\widehat{\underline{\beta}^{rc
}_{1}}
/\{1-0.1(\widehat{\underline{\sigma}}^{2}_{u}+\widehat{\underline{\sigma}}^{2}_{v})\}. \nonumber %\label{mis}
\end{eqnarray}
The rigorous derivation of the relationship is beyond the scope of the paper and we leave it for future research.

Variational approximations to the composite log-likelihood function of crossed random effects models establish a connection between two distinct topics in statistics. We introduce the row-composite likelihood function \eqref{rowc} to eliminate the row dependence of responses by disregarding the column random effects. Subsequently, we utilize the Gaussian variational approximation to eliminate the column dependence of responses through the variational distributions of row random effects. Similarly, we construct the column-composite likelihood function \eqref{columnc},  to eliminate the column dependence of responses by disregarding the random row effects. We then apply the Gaussian variational approximation to eliminate the row dependence of responses through the variational distributions of column random effects. Comparing the row-column composite likelihood-based Gaussian variational approximation and likelihood-based Gaussian variational approximation, both approaches share the common goal of breaking the dependence of responses through manipulating random effects for the generalized linear mixed models with crossed random effects.
For any other links between composite likelihood and variational approximations, we will leave them for future researc
%\section*{Supplementary Material}
%Supplementary material available at \Bka\ online includes proofs of Theorems 1 and 2, additional %simulation results, and additional analyses of the insurance data.

\section*{Acknowledgment}
The authors are grateful to Dr. Fia Fridayanti Adam for sharing the motor vehicle insurance data. This research was partially supported by the Natural Sciences and Engineering Council of Canada.

\section{BibTeX}

\bibliographystyle{Chicago}
\bibliography{paper-ref}
\end{document}